\input amsppt.sty
\def\RefKeys{Fu G1 GG KT LeT Ma Pa Ra Re T1 T2 T3  }
        \widestnumber\key{GG}


\define \OO{\Cal O}

\define \C{\Bbb C}

\define \PP{\Bbb P}

\define \ord{\roman{ ord }}

\define \Bl{\roman{ Bl}}
\define \NB{\roman{ NB}}

\define \Cyc {\roman{ Cyc}}
\define \mult {\roman{mult}}


\define \sub{\subseteq}

\define \inv{^{-1}}
\define \pd{\partial}
\define \pdz#1{\pd_{z_#1}}

\define \pda#1{\pd_{a_{#1}}}




\def\TheMagstep{\magstep1}      
 \def\PaperSize{letter}         

\def\TRUE{TRUE}
\ifx\DoublepageOutput\TRUE \def\TheMagstep{\magstep0} \fi
\mag=\TheMagstep
\TagsOnRight
\vsize48pc
\abovedisplayskip 6pt plus6pt minus0.25pt
\belowdisplayskip=\abovedisplayskip
\abovedisplayshortskip=0mm
\belowdisplayshortskip=2mm
\def\centertext
 {\hoffset=\pgwidth \advance\hoffset-\hsize
  \advance\hoffset-2truein \divide\hoffset by 2
  \voffset=\pgheight \advance\voffset-\vsize
  \advance\voffset-2truein \divide\voffset by 2
 }
\newdimen\pgwidth	\newdimen\pgheight
\def\letter{letter}	
\ifx\PaperSize\letter
	\pgwidth=8.5truein \pgheight=11truein \centertext \fi
\ifx\PaperSize\AFour
	\pgwidth=210truemm \pgheight=297truemm \centertext \fi

\def\atref#1{\ref\key\csname#1\endcsname}
 \newcount\refno \refno=0
 \def\MakeKey{\advance\refno by 1 \expandafter\xdef
 	\csname\TheKey\endcsname{{%
	\ifx\UseNumericalRefKeys\TRUE
		\number\refno
	\else\TheKey \fi}}\NextKey}
 \def\NextKey#1 {\def\TheKey{#1}\ifx\TheKey\NoKey\let\next\relax
  \else\let\next\MakeKey \fi \next}
 \def\NoKey{*!*}
 \expandafter\NextKey \RefKeys *!*
 \ifx\UseNumericalRefKeys\TRUE \widestnumber\no{\number\refno}\fi

 \newskip\sectskipamount \sectskipamount=0pt plus30pt
 \def\sectionhead#1 #2\par{\def\sectno{#1}\def\sectname{#2}%
   \vskip\sectskipamount\penalty-250\vskip-\sectskipamount
   \bigskip
   \centerline{\smc \number\sectno.\enspace\sectname}\nobreak
   \medskip
   \message{\number\sectno. \sectname }%
}

 \def\proclaimm#1#2 {\medbreak\noindent {\bf(\sectno.#2) #1.}\quad
	\begingroup\it}
 \def\subh#1#2{\medbreak\noindent {\bf(\sectno.#2) #1.}\quad}
\def\art#1 #2\par{\subh{\rm({\it#2\unskip\/})}{#1}}

\def\stp{\subh{Setup}}

\def\prop{\proclaimm{Proposition}}
	
\def\thm{\proclaimm{Theorem}}
	
\def\cor{\proclaimm{Corollary}}
\def\lem{\proclaimm{Lemma}}
	
\def\pf{\endgroup\medbreak\noindent {\bf Proof.}\quad}
\def\demobox{\vbox{\hrule\hbox{\vrule\kern.5ex
	\vbox{\kern1.2ex}\vrule}\hrule}}
\def\enddemo{{\unskip\nobreak\hfil\penalty50
  \hskip1em\hbox{}\nobreak\hfil\demobox
  \parfillskip=0pt \finalhyphendemerits=0 \par}}

\def\tgs{\tag\sectno.}

\def\today{\number\day \space\ifcase\month\or
 January\or February\or March\or April\or May\or June\or
 July\or August\or September\or October\or November\or December\fi
\space \number\year}

\magnification=\magstep1
\centertext \topmatter
\title Mixed Segre Numbers and Integral Closure of Ideals
\endtitle

\author Robert Gassler  \endauthor
\abstract
We introduce mixed Segre numbers of ideals which generalize the
notion of mixed multiplicities of ideals of finite colength and show how many
results on mixed multiplicities can be extended to results on mixed Segre
numbers. In particular, we give a necessary and sufficient condition in
terms of these numbers for two ideals to have the same integral closure.
Also, our theory yields a new proof of a generalization of Rees' theorem
that links the integral closure of an ideal to its multiplicity. Finally, we
give a quick application of our results to Whitney equisingularity.
\endabstract

\address  Department of Mathematics, 567 Lake, Northeastern University,
Boston, MA 02115, U.S.A., e-mail: gassler\@neu.edu  \endaddress
\thanks It is a pleasure to thank Terry Gaffney for many helpful
conversations. I also thank him for encouraging me to improve the results of
an earlier version of this work.\endthanks
\endtopmatter

\document
\define\intro{0}
\define \mixed{1}
\define \segre{2}
\define \surf{3}
\define \mink{4}
\centerline{alg-geom/9703018}\medskip

\define\th#1{\tilde h_{#1}}
\define\tD#1{\tilde D_{#1}}
\define\tb#1{\tilde b_{#1}}


\sectionhead {\intro} Introduction


In 1973, Bernard Teissier \cite{T1} used the multiplicity of ideals that are
primary to the maximal ideal $m$ of a local ring to study the geometry of
isolated hypersurface singularities. In fact, to study Whitney conditions,
one has to understand the integral closure of the product of the Jacobian
ideal and the maximal ideal. Now, Rees showed in his celebrated theorem that
two $m$--primary ideals $I\sub J$ of a formally equidimensional local
ring have the same integral closures if and only if their
multiplicities are equal; see e.g. \cite{T1}. Teissier gave a formula that
expresses the multiplicity of the product of $m$ and an $m$--primary ideal in
terms of some mixed multiplicities.

Later, in \cite{T2}, Teissier defined mixed multiplicities for any pair of
$m$--primary ideals. And, in \cite{T3}, mixed multiplicities were
used to give a numerical criterion for two $m$--primary ideals to have the
same integral closure.

\medskip
More recently, Terry Gaffney and the author \cite{GG} defined Segre numbers of
ideals to extend Teissier's result on isolated hypersurface singularities to
non--isolated ones. Segre numbers are a natural generalization of
multiplicities of $m$--primary ideals. In fact, they allowed us to extend
Rees' theorem to arbitrary ideals. It seems natural to define also mixed
Segre numbers for any pair of ideals, and to use them to find
generalizations of Teissier's results. Indeed, Corollary (\mink.4) of this
work gives a necessary and sufficient criterion in terms of mixed Segre
numbers for the integral closure of two ideals to coincide.

\medskip
Teissier also proved inequalities relating the multiplicity of the
product of two ideals and mixed multiplicities. To prove this, Teissier
showed that it is enough to prove such an inequality on a normal surface, and
used then the negative definiteness of the intersection matrix of the
components of the exceptional divisor of a resolution of singularities of the
surface.

\medskip
Following Teissier's approach, we show that similar inequalities hold for
Segre numbers. However, for such an inequality of Segre numbers of
codimension $k$ to hold, the ideals in question need to have the same
behavior in lower codimension. This, in turn, can be ensured by a condition
on the mixed Segre numbers of codimension less than $k.$ Again, the proof
uses certain properties of the resolution of singularities of normal surfaces
(see Section {\surf}). Interestingly, the behaviour of the integral closure
of an ideal in codimension $k$ influences its behaviour in higher
codimension; see the remark after Corollary (\surf.4).

\medskip
Finally, Teissier used his product formula to derive Minkowski--type
inequalities that relate the multiplicity of the product of the two ideals
with their individual multiplicities. Similar product formulas are proven
here in Section {\mink}. Also, we extend Teissier's result that these
inequalities are equalities if and only if some powers of the two ideals have
the same integral closure. However, for general ideals, we need to assume
that the Segre cycles of one ideal satisfy a chain condition; see Theorem
(\mink.8).

\medskip
Finally, we give a new proof of the recent generalization \cite{GG, Corollary
(4.9)} of Rees' theorem.

\medskip
This paper shows that Segre numbers are a good generalization of multiplicity
of ideals. Also, by using `induction on the codimension,' most of the results
for ideals of finite colength can be generalized without too many problems to
ideals of lower codimension.

\medskip
The paper is organized as follows. In Section {\mixed}, we review Teissier's
theory of mixed multiplicities. Using an easy combinatorial result
(\mixed.2), we generalize his results to arbitrary local rings of
equidimensional analytic space germs.

Mixed Segre numbers are introduced in
Section {\segre}. Then, in Section {\surf} we study the case of ideals on
surfaces. This is the model case which is used in Section {\mink} to study the
general case. For this, we first review an alternative definition of Segre
numbers and cycles, and then examine moving and fixed components of these
cycles in more detail in (\mink.1) and (\mink.2). The criterion for two ideals
to have the same integral closure is then given in Corollary (\mink.4). Then,
we discuss product formulas and Minkowski--type inequalities. Finally, a first
application to equisingularity theory is given in (\mink.11).


\sectionhead {\mixed} Mixed Multiplicities


\stp1 Let $(X,0)\sub(\C^N,0)$ be an analytic germ of pure dimension $n.$ We
denote the maximal ideal of its local ring $\OO_{X,0}$ at 0 by $m.$  Consider
two $m$--primary ideals $I_1,I_2$ in this local ring. Then the multiplicity
$e(I_1)$ of
$I_1$ equals the dimension of the $\C$--vector space
$$\OO_{X,0}/(f_1,\dots,f_n)$$ where $f_1,\dots,f_n$ are generic linear
combinations of generators of $I_1.$ Inspired by this theorem of Samuel,
Teissier defined mixed multiplicities of two ideals as follows:
$$e_{i,n-i}(I_1,I_2)=\dim\OO_{X,0}/(f_1,\dots,f_i,g_1\dots,g_{i-n}),$$
where the $f_j$ are generic linear combinations of generators of $I_1$ and
the $g_k$ are generic linear combinations of generators of $I_2.$ Clearly,
this definition implies $e_{i,n-i}(I_1,I_2)=e_{n-i,i}(I_2,I_1).$

In \cite{T2}, Teissier showed that for normal $X$ the following inequality
holds.
$$e_{i,n-i}(I_1,I_2)^n\leq e(I_1)^ie(I_2)^{n-i}.\tgs1.2$$
Furthermore, in \cite{T3}, he gave a criterion in terms of mixed
multiplicities that determines when the integral closure of two $m$--primary
ideals coincide. To prove a slight generalization of this result, we need the
following observation.

\lem2 Let $(a_1,\dots,a_k),(b_1,\dots,b_k),(c_1,\dots,c_k)$
$k$--tuples of nonnegative integers. Denote the sum of their elements by
$a,b,c,$ respectively. Suppose that for
$i=1,\dots,k$ the inequality $a_i^2\leq b_ic_i$ obtains. Then $a=b=c$ if, and
only if, for $i=1,\dots,k$ we have $a_i=b_i=c_i.$

\pf Consider the following inequalities.
$$a=\sum_{i=1}^k a_i\leq \sum_{i=1}^k\sqrt{b_ic_i}\leq \sqrt{b}
\sqrt{c}.$$ The second inequality is the Cauchy--Schwarz
inequality. It is an equality if, and only if, the tuples $b$ and $c$ are
linearly dependent. In our situation, this is equivalent to $b_i=c_i.$ The
claim follows.
\enddemo

\medskip
We shall prove now Teissier's result \cite{T3, p.354} in a slightly more
general situation.

\thm3 Let $(X,0)$ be an equidimensional complex analytic space germ of
dimension $n.$ Let $I_1,I_2$ be two ideals in the local ring
$\OO_{X,0},$ primary with respect to the maximal ideal $m.$ Then their
integral closures coincide if, and only if, we have
$$e(I_1)=e_{n-1,1}(I_1,I_2)=\dots=e_{1,n-1}(I_1,I_2)=e(I_2).$$

\pf Teissier showed that it is enough to prove the result for 2--dimensional
germs. (In his reduction argument \cite{T3, p.348} he assumes that $X$ is
Cohen-Macaulay. It is easy to see  that this assumption is unnecessary.) In
his original statement, Teissier also assumes that $X$ is normal, and
proves then the statement for normal surface germs using a result
of Ramanujam \cite{P} that expresses mixed multiplicities in terms of orders
of vanishing along components of the exceptional set of a resolution of
singularities of the surface germ.

The above Lemma will allow us to conclude the statement for non--normal
surfaces  from Teissier's results. This is the only new ingredient.

We assume now that $X$ is a surface germ. Let $n:\bar X\to X$ be its
normalization, and assume that the preimage of 0 is formed by points
$x_1,\dots, x_k.$ Consider the ideals $I_1^{(i)}$ and $I_2^{(i)}$ in
the local rings $\OO_{\bar X,x_i}$ induced by the pullbacks of the ideals on
$X.$ Then, any path $\phi:(\C,0)\to(X,0)$ lifts to a path
$\bar\phi:(\C,0)\to(\bar X,0).$ Therefore, by the valuative criterion for
integral closure, we have that $I_1$ and $I_2$ have the same integral closure
if, and only if, for $i=1\dots,k$ the integral closure of $I_1^{(i)}$ and
$I_2^{(i)}$ coincides. By Teissier's result, this happens if, and only if, we
have
$$e(I_1^{(i)})=e_{1,1}(I_1^{(i)},I_2^{(i)})=e(I_1^{(i)})$$ for all $i.$
Using a projection formula for multiplicities \cite{F, (4.3.6)} and
(\mixed.1.1), we see
$$e(I_1)=\sum_{i=1}^ke(I_1^{(i)}),\quad e(I_2)=\sum_{i=1}^ke(I_2^{(i)}),
\quad e_{1,1}(I_1,I_2)=\sum_{i=1}^ke_{1,1}(I_1^{(i)},I_2^{(i)}).$$
Our assumption, together with Lemma (\mixed.2), imply that we can apply
Teissier's result to the normal surface $\bar X$ and the ideals induced by
$I_1$ and $I_2.$ This finishes the proof.
\enddemo

\cor4 In the Setup (\mixed.1), the following inequalities obtain for
$i=1,\dots,n.$
$$\align
e_{i,n-i}(I_1,I_2)^n&\leq e(I_1)^ie(I_2)^{n-i},\tgs4.1\\
e(I_1I_2)^{1/n}&\leq
e(I_1)^{1/n}+e(I_2)^{1/n}.\tgs4.2\endalign$$ The second inequality is an
equality if, and only if, there exist postive integers $a,b$ so that the
integral closures of
$I_1^a$ and $I_2^b$ coincide.

\pf The first inequality follows, using Teissier's arguments \cite{T2}, from
the the following inequality in the local ring of a surface singularity:
$$e_{1,1}(I_1,I_2)^2\leq e(I_1)e(I_2).$$ Now, with the notation of the above
proof, we have
$$\eqalign{e_{1,1}(I_1,I_2)^2=&\left(\sum e_{1,1}(I_1^{(i)},I_2^{(i)})\right)^2
\leq\left(\sum\sqrt{e(I_1^{(i)})e(I_2^{(i)})}\right)^2\cr
\leq&\left(\sum e(I_1^{(i)})\right)\left(\sum e(I_2^{(i)})\right).}$$ Here,
the first inequality is Teissier's result, the second is the Cauchy--Schwarz
inequality.

Finally, the second inequality (\mixed.4.2) follows now from Teissier's
expansion formula for multiplicities. Furthermore, Teissier's proof shows the
last assertion if we replace his argument
\cite{T3, p. 354} by the above theorem.\enddemo


\sectionhead {\segre} Mixed Segre Numbers


\stp1 Let $(X,0)\sub(\C^N,0)$ be an equidimensional analytic germ of pure
dimension $n.$ Recall from \cite{GG, (2.1) and (2.2)} that  polar varieties
$P_k(I,X)$ and Segre cycles $\Lambda_{k+1}(I,X)$ of an ideal
$I$ of nowhere dense co--support on $X$ are defined inductively for
$k=0,\dots,n-1$ as follows:
$P_0(I,X)=X$; for $k\geq 1$ we define $P_k(I,X)$ to be the
closure of $V(f|_{P_{k-1}(I,X)})-V(I),$ where $f$ is a general linear
combination of generators of $I.$ The $k$th Segre cycle $\Lambda_k(I,X)$ can
be defined as the difference of cycles
$$[V(f|_{P_{k-1}(I,X)})]-[P_k(I,X)].$$ Although in general both, polar
varieties and Segre cycles, depend on the choice of the general linear
combinations of generators of $I,$ their multiplicities at 0 are
well--defined and denoted by $m_k(I,X)$ and $e_k(I,X).$ We will often omitt
the $X$ in these notations, and write, for example, $m_k(I)$ for $m_k(I,X).$

For a subspace $Y$ of $X$ no
component of which is contained in $V(I),$ we define the $k$th Segre number
$e_k(I,Y)$ of $I$ on $Y$ to be $e_k(I',Y),$ where $I'$ is the ideal induced
by $I$ in the local ring of $Y$ at 0.

We refer to \cite{GG} for the general development of the theory of Segre
numbers and cycles.

\medskip
In this work, we will mainly use the following two properties of Segre
numbers:

(1) They are well--behaved with respect to polar varieties
(\cite{GG,(2.2.1)}), i.e. for $i+j\leq n$ and $i>0,$ we have
$$e_{i+j}(I)=e_i(I,P_j(I)).$$

(2) The one--codimensional Segre number can be computed on a plane section
off 0 (\cite{GG,(2.3)}): Let
$p:\C^N\to\C^{n-1}$ be a general linear projection, $\epsilon$ a general
point of $\C^{n-1}$ close
to 0, and
$H=p\inv(\epsilon).$ Then, we have
$$e_1(I)=\sum_{x\in V(I)\cap H}e(I,(X\cap H,x)).$$

Now, let
$I_1,I_2$ be ideals in the local ring $\OO_{X,0}$ of codimension at least
one. In the following, we fix
$n$--tuples $\bold f=(f_1,\dots,f_n)$ and $\bold g=(g_1,\dots,g_n)$ of linear
combinations of generators of $I_1$ and $I_2$ respectively. For $k=1,\dots,n$
we will denote the $k$--tuples $(f_1,\dots,f_k)$ and $(g_1,\dots,g_k)$ by
$\bold f^{[k]}$ and $\bold g^{[k]}.$ We denote the ideal generated by $\bold
f^{[k]}$ by $I_1^{[k]},$ and define $I_2^{[k]}$ analogously.

We will first examine how many generic generators of an ideal are necessary to
define polar varieties and Segre cycles.

\prop2 In the above setup, outside the $k$--codimensional polar variety
$P_k^{\bold f}(I_1)$ of $I_1,$ the ideal $I_1^{[k]}$ is a reduction of
$I_1.$

\pf Suppose $I_1$ is generated by $m+1$ elements, and consider the blowups
$B=\Bl_{I_1}X\sub X\times \PP^m$ and $B'=\Bl_{I_1^{[k]}}X\sub X\times
\PP^{k-1}$ Then, the identity map on $X$ induces a map
$$B-(B\cap (X\times K)\to B'$$ induced by a central projection
$\pi:\PP^m-K \to \PP^{k-1}$ with $(m-k)$--dimensional center $K.$ Now, as
$\bold f$ is formed by generic linear combinations, the center $K$ is is
general position with respect to $X.$ So, a standard argument shows that
$B\cap (X\times K)$ is of minimal dimension $n-k.$ The image of this
intersection in
$X$ equals
$P_k^{\bold f}({I_1}).$

Note that if $x$ is point of $X$ where the fibre $B(x)$ of $B$ over $x$ is of
dimension bigger than $k,$ then $B(x)\cap K$ is non--empty. Hence, the point
$x$ is contained in the mentioned polar variety.

Also, if $x$ is outside the polar variety, the induced projection $B(x)\to
\PP^{k-1}$ is finite. Otherwise, there would be a curve $C$ in $B(x)$ that is
mapped to one point $l$ in $\PP^{k-1}.$ Now, $\pi\inv(l)$ is isomorphic to
$\C^{m-k+1}=\PP^{m-k+1}-K.$ As $B(x)$ is closed, the closure $\bar
C\sub\PP^{m-k+1}$ is also contained in $B(x),$ and, by Bezout's theorem,
intersects $K.$ This is a contradiction to the assumption on $x.$

Therefore, the induced projection
$\Bl_{I_1}(X-P_k^{\bold f}({I_1}))\to
\Bl_{I_1^{[k]}}(X-P_k^{\bold f}({I_1}))$ is finite and generically
one--to--one. So, by a well--known characterization of integral dependence
(see e.g.
\cite{T1, Ch.0, 0.4, p.288}), the result follows.\enddemo

\cor3 In the above setup, we have the equalities:
$$P_k^{\bold f}(I_1)=P_k^{\bold f^{[k+1]}}(I_1^{[k+1]}),\qquad
\Lambda_k^{\bold f}(I_1)=\Lambda_k^{\bold f^{[k+1]}}(I_1^{[k+1]}).$$

\pf The integral closures of $I_1$ and $I_1^{[k+1]}$ coincide outside
$P_{k+1}^{\bold f}(I_1),$ therefore the desired equalities hold on
$X-P_{k+1}^{\bold f}(I_1).$ Now, the $(k+1)$--codimensional polar variety,
which was cut out, is nowhere dense in the $k$--dimensional objects in
question, so the equality remains when passing to the closure.\enddemo

\medskip
Looking back into the definition, we see that the element $f_{k+1}$ is used
in the formation of $\Lambda_k(I)$ only to distinguish components of $P_k(I)$
and
$\Lambda_k(I)$ in the cycle $V(f_k|P_{k-1}(I)).$ So, if we know the support
of $V(I_1)\cap P_{k-1}(I),$ then $f_{k+1}$ is not needed to define the Segre
cycle of codimension $k.$ This observation will lead us to the correct
definition of mixed Segre cycles.

\art4 Mixed Segre numbers

Before we define
mixed Segre numbers in general, we consider the special case of
2--codimensional mixed Segre cycles and numbers. The first naive approach,
following Teissier's definition in the case of ideals of finite colength,
would be to consider the Segre cycle of codimension 2 of
$I_1^{[1]}+I_2^{[1]}.$ However, as we have seen in Corollary (\segre.2), we
generally need three elements to define the Segre cycles of codimension 2,
unless we give ourselves some set that will play the role of the support
of $I$ as in the above observation. So, we define the mixed Segre cycle as
folllows: Let $h_1$ and $h_2$ be two generic linear combinations of $f_1$ and
$g_1.$ Define $P_1^{1,1}(I_1,I_2)$ to be the closure of $V(h_1)-V(I_1+I_2),$
and $P_2^{1,1}(I_1,I_2)$ to be the closure of
$V(h_2|P_1^{1,1}(I_1,I_2))-V(I_1+I_2).$ Finally, define the Segre cycle
$\Lambda_2^{1,1}(I_1,I_2)$ to be the part of the 2--codimensional cycle
$[V(h_2|P_1^{1,1}(I_1,I_2))]$ supported in $V(I_1+I_2).$

In general, for $k=1,\dots,n,$ and positive integers
$i,j,$ we define the mixed Segre number of codimension $k$ as follows: Let
$\bold h$ be a
$k$--tuple of generic linear combinations of the elements
$f_1,\dots,f_i,g_1,\dots,g_j.$ Then, we define inductively $P_k^{i,j}(I_1,I_2)$
to be the closure of $V(h_k|P_{k-1}^{i,j}(I_1,I_2))-V(I_1+I_2)$
and $\Lambda_k^{i,j}(I_1,I_2)$ to be the part of the cycle
$[V(h_k|P_{k-1}^{i,j}(I_1,I_2))]$ supported in $V(I_1+I_2).$ Finally, let
$e_k^{i,j}(I_1,I_2)$ be its multiplicity.

We also define $e_k^{k,0}(I_1,I_2):=e_k(I_1)$ and
$e_k^{0,k}(I_1,I_2):=e_k(I_2).$

In particular, for ideals of finite colength our definition of mixed
multiplicities coincide with Teissier's.


\sectionhead {\surf} The Surface Case


We will first discuss ideals on a normal surface $X.$
We will also use classical results on resolutions of normal surfaces, which we
are going to review.

\art1 Resolutions of normal surfaces

An outline of the proofs and more
references can be found in Fulton's book \cite{F,  Ex. 2.4.4 , p.39 and Ex.
7.1.16, p.125}.

Let
$(X,0)$ be a normal surface, and
$\pi:\tilde X\to (X,0)$ a resolution of singularities. Denote the irreducible
components of the exceptional fibre $\pi\inv(0)$ by $E_1,\dots,E_r.$ Then, for
an irreducible curve
$C\subset (X,0)$, there exist unique positive $a_j\in\Bbb Q$ so that for all
$i=1,\dots,r$ the equality of zero--cycles in $\pi\inv(0)$ obtains:
$$\tilde C\cdot E_i+\sum_{j=1}^ra_j(E_j\cdot E_i)=0.$$ Here, $\tilde C$ is
the strict transform of $C$ by $\pi.$ We define
$$C'=\tilde C+\sum_{j=1}^ra_jE_j.$$ Then, if
$D$ is a divisor on
$X,$ we have $[D']=[\pi^*D].$

The main ingredient for the proof is the negative definiteness of the
intersection matrix $\{(E_i,E_j)\}_{i,j},$ where $(E_i,E_j)$ denotes the
degree of the intersection product $E_i\cdot E_j.$ This degree is
well--defined as the intersection is supported in the complete fibre
$\pi\inv(0)$ of
$\tilde X$ over 0.

We will extend the notion of strict transforms to cycles on $X$ by linearity;
again the strict transform of a cycle $S$ by $\pi$ will be denoted by $\tilde
S.$

\medskip For the proof of inequality of mixed Segre numbers on a surface
germ, we will need the following modified version of the Cauchy--Schwarz
inequality.

\lem2 Let $u,v,w$ be elements of a real vector space $V,$ and assume
$\langle,\rangle$ is a positive definite bilinear form on $V.$ Then, if
$\langle u,w\rangle\geq\langle v,w\rangle \geq 0$ holds, we have the
following inequality
$$\langle u+w,v\rangle ^2\leq \langle u+w,u\rangle \langle v+w,v\rangle .$$

\pf Let $a=u-v,$ and consider the function
$$f(t)=\langle v+ta+w,v+ta\rangle \langle v+w,v\rangle -\langle
v+ta+w,v\rangle ^2.$$ Clearly, $f(0)=0;$ we are going to show that its
derivative $f'(t)$ is positive for $t\geq 0.$ This implies the desired
inequality.

Now, an easy computation shows
$$\eqalign{f'(t)=2t(\langle a,a\rangle \langle w,v\rangle &+\langle
a,a\rangle \langle v,v\rangle -\langle a,v\rangle ^2)\cr
&+\langle a,w\rangle \langle v+w,v\rangle .}$$ Now, the
Cauchy--Schwarz inequality together with the assumptions implies $f'(t)\geq
0.$ This finishes the proof.\enddemo

\prop3 Assume that $X$ is a normal surface and the equality of
one--cycles
$$[V(I_1)]_1=[V(I_2)]_1$$ hold, where $[Y]_1$ denotes the one--cycle formed by
the one-dimensional components of the cycle $[Y].$ Then, the following
inequality
obtains:
$$e_2^{1,1}(I_1,I_2)^2\leq e_2(I_1)e_2(I_2).$$

\pf First, we claim the inequality $e_2^{1,1}(I_1,I_2)\leq
e(I_2,P_1(I_1)).$ In fact, consider a generic linear combination
$h=f_1+tg_1$ of the generators of $I_1^{[1]}+I_2^{[1]}.$ Then, we
have
$$I_1^{[1]}+I_2^{[1]}|P_1^h(I_1^{[1]}+I_2^{[1]})=
I_2^{[1]}|P_1^h(I_1^{[1]}+I_2^{[1]}),$$
as $f_1|P_1^h(I_1^{[1]}+I_2^{[1]})=-tg_1|P_1^h(I_1^{[1]}+I_2^{[1]}).$ Next,
consider the space
$$Y=\overline{V(f_1+tg_1)-V(I_1)\times\C}\subset X\times\C$$ where $t$ is now
the coordinate on $\C.$ Then, as we have the equality of cycles
$$[V(I_1)]=[V(I_2)]=[V(I_1^{[1]}+I_2^{[1]})],$$ the fibre of $Y$ over $t=0$
equals $P_1(I_1),$ while the fibre over non--zero $t$ equals
$P_1^{f_1+tg_1}(I_1^{[1]}+I_2^{[1]});$ see \cite{GG, (4.3)} for a proof of
this. So, by the upper semicontinuity of the multiplicity, and as $g_1$ was
chosen generically, we have
$$e_2^{1,1}(I_1,I_2)=e(I_2^{[1]},P_1^h(I_1^{[1]}+I_2^{[1]}))\leq
e(I_2^{[1]},P_1(I_1))=e(I_2,P_1(I_1)).$$

Next, consider a resolution of singularities
$\pi:\tilde X\to X$ as above with the additional property that $\pi^*I_1$ and
$\pi^*I_2$ are invertible sheaves, e.g. a resolution of singularities of the
blowup of
$X$ along
$I_1I_2.$ Note that by the normality of $X,$ this resolution is still an
isomorphism outside the exceptional fibre. In fact, as $X$ is normal, it is
smooth outside 0, and so at a point $p$ outside 0 the ideal induced by $I_1$
in $\OO_{X,p}$ is principal. Hence, the blowup of $(X,p)$ along $I_1$ is
isomophic to $(X,p).$ Also, we assume that the
strict transforms of $P_1(I_1)$ and
$\Lambda_1(I_1)$ in $\tilde X$ don't meet. This can be achieved by replacing
$\tilde X$ by a successive blowup of finitely many points in $\tilde X.$

Then, by the above, there exist positive numbers
$u_i,v_i,w_i\in\Bbb Q$ so that we have for all $j$
$$\align
\left(\Lambda_1(I_1)\tilde{\,}+\sum w_iE_i,E_j\right)&=0,\\
\left(P_1(I_1)\tilde{\,}+\sum u_iE_i,E_j\right)&=0\\
\left(P_1(I_2)\tilde{\,}+\sum v_iE_i,E_j\right)&=0.\endalign$$
Furthermore, consider the divisor on $X$ given by $f_1=0$ with associated
Weil divisor $[P_1(I_1)]+\Lambda_1(I_1).$ Then, by the above we have
$u_i+w_i=\ord_{E_i}\pi^*f_1.$

Now, by the projection
formula and the assumed properties of $\tilde X,$ we have
$$\eqalign{
e_2(I_1)&=e(\pi^*I_1,P_1(I_1)\tilde{\,})=
(\Lambda_1(I_1)\tilde{\,}+\sum(\ord_{E_i}\pi^*(I_1)E_i,P_1(I_1)\tilde{\,})\cr
&=(\sum(\ord_{E_i}\pi^*I_1)E_i,P_1(I_1)\tilde{\,})
=-\sum_{i,j}(u_i+w_i)u_j(E_i,E_j).}$$ For the last equality, we used that the
pullback of $I_1$ is invertible; hence its order of vanishing along $E_i$
equals the order of vanishing of $f_1$ along $E_i.$

Similarly, we obtain
$$\align
e_2(I_2)&=-\sum_{i,j}(v_i+w_i)v_j(E_i,E_j)\text{ and}\\
e(I_2,P_1(I_1))&=-\sum\sum_{i,j}(v_i+w_i)u_j(E_i,E_j).
\endalign$$
Hence, it is enough to show the inequality $\langle v+w,u\rangle ^2\leq\langle
u+w,u\rangle \langle v+w,v\rangle ,$ where $\langle\, ,\rangle $ denotes the
positive definite bilinear form given by the matrix
$\{-(E_i,E_j)\}_{i,j}.$

After exchanging $I_1$ and $I_2$, we may assume $\langle u,w\rangle
\geq \langle v,w\rangle.$ Hence, to apply the above lemma, it remains to check
$\langle v,w\rangle \geq 0.$ But, we have
$$\langle v,w\rangle =(\sum
v_iE_i,\Lambda(I_1)\tilde{\,})> 0,$$
as $(\Lambda_1(I_1)\tilde{\,},E_i)\geq 0$ for all $i.$ Furthermore, as
$\Lambda_1(I_1)$ is a cycle through the origin, its strict transform
intersects at least one component of the exceptional fibre. Hence, for at least
one index
$i,$ the intersection number is non--zero. This finishes the proof.
\enddemo

\cor4 In the situation of the above Proposition (\sectno.3), the integral
closures of $I_1$ and $I_2$ coincide if and only if the equalities
$$e_2(I_1)=e_2^{1,1}(I_1,I_2)=e_2(I_2)$$ obtain.

\pf If the two ideals have the same integral closure, then the equalities of
2--codimensional Segre numbers clearly obtain. So, assume now that the
equalities hold. If both, $I_1$ and $I_2,$ are of finite colength, this is
a result of Teissier \cite{T3, p.354}. So, we consider now the case where
$I_1$ and $I_2$ are of codimension one; in other words, their
one--dimensional Segre cycles are not 0. We also assume that at least one of
the two--codimensional Segre numbers in question is not zero. Otherwise, both
ideals are principal, and then, by the assumptions, obviously equal.

The equalities imply that the inequality in Proposition (\sectno.3) is an
equality. Consider now the derivative $f'(t)$ of the function $f(t)$ in the
proof of Lemma (\sectno.2) that was used to prove this equality. It has to
be identically to
0 as, by the assumptions, $f(0)=f(1)$ and $f'(t)\geq 0.$ As we have seen in the
above proof, the product
$\langle w,v\rangle$ is positive. Also, by the classical
Cauchy--Schwarz inequality,
$\langle a,a\rangle
\langle v,v\rangle-\langle a,v\rangle^2$ is positive. Hence $\langle
a,a\rangle=0,$
and so $a=0.$ Hence, by a characterization of integral closure
\cite{LT}, the two ideals have the same integral closure. This finishes the
proof.\enddemo

\medskip
Note that, for ideals of codimension one, in order to prove that the equalities
imply that the ideals have the same integral closure, we only used the fact
that the assumed equalities imply that the inequality in Proposition
(\sectno.3) is actually an equality. This should be contrasted to Teissier's
results for ideals of finite colength: In this case, if the inequality is an
equality, we can only conclude that some powers of the ideals have the same
integral closure.

The stronger result for ideals of codimension one comes from the assumption
that we assumed the ideals to have the same integral closure outside 0. As the
proof of the inequality of Proposition (\sectno.3) showed, the behavior
outside 0
influences their behavior at 0. We will use this observation later on to derive
some results on Minkowski--type inequalities.

\medskip
From the above proposition and the corollary, we can also derive similar
results for arbitrary surfaces.

\cor5 Let $(X,0)$ be an arbitrary surface germ, and $I_1,J$ be two ideals of
codimension at least one in the local ring of $X$ at 0. Assume that the
integral closure of $I_1$ and $I_2$ at points of $X$ outside 0 coincide
(that is for
suitably chosen representatives.) Then the inequality
$$e_2^{1,1}(I_1,I_2)^2\leq e_2(I_1)e_2(I_2)$$ obtains. Furthermore, $\bar
I_1=\bar I_2$ if and only if we have the equality
$$e_2(I_1)=e_2^{1,1}(I_1,I_2)=e_2(I_2).$$

\pf Again, by (\mixed.3), it is enough to prove this
ideals of codimension one. Consider the normalization
$n:(\bar X,S)\to (X,0)$ where
$S$ is a finite set of points. Then, the integral closures of $I_1$ and $I_2$
coincide if and only if the integral closures of the pullbacks $n^*I_1$ and
$n^*I_2$ coincide. Furthermore, by the above Proposition (\surf.3), at each
point $x$ of $S$ we have an inequality
$$e_2^{1,1}(n^*I_1I,n^*I_2,(\bar X,x))^2\leq e_2(n^*I_1,(\bar
X,x))e_2(n^*I_2,(\bar X,x)).$$ Hence, using Cauchy--Schwarz inequality,
we get the desired inequality. Also, we see, using Lemma (\mixed.2), that the
equality in the statement of the corollary obtains if and only if the above
equality holds at each point
$x$ of
$S.$ So, the claim follows from the above proposition.\enddemo

Finally, we can also give a numerical criterion for the
integral closure of
$I$ and
$I_2$ to coincide outside 0:

\cor6 Let $(X,0)$ be an arbitrary surface germ, and $I_1,I_2$ be two ideals of
codimension at least one in the local ring of $X$ at 0. Then $\bar I_1=\bar
I_2$ if and only if the following equalities obtain:
$$e_2(I_1)=e_2^{1,1}(I_1,I_2)=e_2(I_2)\text{ and
}e_1(I_1)=e_1^{1,1}(I_1,I_2)=e_1(I_2)$$

\pf We are going to show that the second equality implies that $I_1$ and $I_2$
have the same integral closure outside 0. By symmetry, it is enough to show
that $I_2$ is integrally dependent on $I_1$ outside 0. So, consider the
normalized blowup $b:\NB_{I_1}X\to X$ with exceptional divisor $D.$ Then, we
need to show that for a component $C$ of $D$ that doesn't map to 0, the order
of vanishing of $b^*I_2$ along $C$ is no smaller than the order of $D$ along
$C.$ This is equivalent to the following claim: For a generic hyperplane $H$
off 0 that intersects
$b(C)$ transversally, the ideal $I_2|X\cap H$ is integrally dependent on
$I_1|X\cap H.$ In fact, as $H$ is chosen generically, its preimage $b\inv(H)$
intersects $C$ transversally, and the order of $b^*I_2$ along $C$ equals the
order of $b^*(I_2|H)$ along $b\inv(H)\cap C.$ The latter is a component of the
exceptional divisor of the blowup of $X\cap H$ along $I_1|X\cap H.$

Now, $(X\cap H,V(I_1)\cap H)$ is a multi-curve germ, and so, by Rees Theorem,
it is enough to show that at each point of $V(I_1)\cap H$ the multiplicities of
the ideals induced by $I_1$ and
$I_1+I_2$ in the local ring of $X\cap H$
coincide. Furthermore, again as we're working on a curve, both
$I_1$ and $I_2$ have a reduction generated by one element. So,
$I_1^{[1]}+I_2^{[1]}$ is a reduction of $I_1+I_2.$ Now, at a point
$x\in V(I_1+I_2)\cap H\sub V(I_1)\cap H,$ we have
$$e(I_1+I_2,(X\cap H,x))\leq e(I_1,(X\cap H,x)).$$
So, as
$e_1(I_1)=e(I_1,X\cap H)$ and similarly for $I_1+I_2,$ the assumptions imply
the equality of sets $V(I_1)\cap H=V(I_1+I_2)\cap H,$ and, furthermore, at
each point of this set the above inequality is actually an equality.
This finishes the proof.\enddemo


\sectionhead {\mink} The general case


We now return to the setup (\segre.1) to study mixed Segre numbers in the
general case. We do this by reducing it again to the the surface case;
however, we will have to deal now with non--local phenomena. The first lemma
is a key--ingredient.

For the following, we need an alternative description of the Segre cycles,
discussed in \cite{GG, (2.1)}. We review this briefly and discuss the
phenomena of moving components in more detail.

\art1 Fixed and moving parts of Segre cycles

In the Setup (2.1), consider the blowup
$$\Bl_IX\subset X\times\PP^m,$$ where $m+1$ is the number of a set of
generators of $I,$ with exceptional divisor $D.$ Denote the part of $D$
formed by components mapping to 0 by $D^0,$ and the part formed by the other
components by $D^{X-0}.$ A hyperplane $H\subset \PP^m$ gives rise, via
pullback by the map $\Bl_IX\to\PP^m$ induced by the projection onto
the second factor, to a Cartier divisor on the blowup, which we denote again
by $H.$ Intersecting with $H$ represents the first Chern class of the
canonical line bundle of the blowup.

Then, let $H_1,\dots,H_{n-1}$ be generic
hyperplanes in the projective space. Consider the cycle
$$b_*(H_1\cap\dots\cap H_k\cap D^{X-0})$$ in $X,$ where $b$ is the canonical
map of the blowup. It is a Segre cycle of $I$ of codimension $k+1.$ Finally,
the degree of the 0--cycle $$H_1\cap\dots\cap H_{n-1}\cap D^0$$ equals the
top Segre number $e_n(I).$

An easy dimension count shows that for a component $C$ of $D,$ the cycle
$$b_*(H_1\cap\dots\cap H_k\cap C)$$ is non--zero if and only if either $b(C)$
is of codimension $k+1$ in $X,$ or the image $b(C)$ is of lower codimension and
the fibre of $C$ over 0 has dimension bigger than $k.$ In the second case,
the non--zero cycle associated to $C$ is called a {\it moving component} of the
Segre cyce $\Lambda_{k+1}(I).$ Indeed, it depends on the choice of the
hyperplanes $H_i.$

Now, consider a generic linear projection $p:\C^N\to\C^{n-k-1}$ and a general
point
$\epsilon$ in its target close to 0. We will call the plane
$L=p\inv(\epsilon)$ of codimension $n-k-1$ a {\it generic
$(n-k-1)$--codimensional plane off 0.} We consider a small
representative of $X,$ again denoted by $X,$ the space $X\cap L,$
and the ideal sheaf induced by a representative of $I.$ If all these
representatives are chosen small enough, the plane $L$ intersects
$\Lambda_{k+1}(I)$ transversally, and the degree of the intersection equals
$e_{k+1}(I).$

As $L$ was chosen generically, its preimage in $\Bl_IX$ also intersects the
exceptional divisor transversally. In fact, $b\inv(L)$ is
isomorphic to the blowup of $X\cap L$ along the ideal induced by $I,$ and
$D\cap L$ is the exceptional divisor of this blowup. So, if $C$ is a
component of $D$ that maps to a subset of codimension $k+1,$ it induced the
component $C\cap L$ of the exceptional divisor of the blowup of $X\cap L.$
Also, as the intersection is transversal and the canonical line bundle of
$\Bl_IX$ restricts to the one of the blowup of $X\cap L,$ the degree of the
0--cycle
$$H_1\cap\dots\cap H_k\cap C\cap b\inv(L)$$
equals the multiplicity of the cycle
$$b_*(H_1\cap\dots\cap H_k\cap C)$$
of codimension $k+1$ in $X.$ If we now sum over all such components, we get
$$\sum_{x\in X\cap L}e_{k+1}(I,(X\cap L,x))=e_{k+1}(I)_f,\tgs1.1$$
where $e_{k+1}(I)_f$ is the contribution to $e_{k+1}(I)$ coming from the
fixed components of the Segre cycle $\Lambda_{k+1}(I);$ we call it the
{\it fixed part} of $e_{k+1}(I).$ Note that the summand on the left--hand side
is non--zero only at finitely many points.

We denote the contribution to $e_{k+1}(I)$ coming from the
moving components of the Segre cycle $\Lambda_{k+1}(I)$ by $e_{k+1}(I)_m,$ the
{\it moving part} of $e_{k+1}(I).$ Then, we have
$$e_{k+1}(I)=e_{k+1}(I)_m+\sum_{x\in X\cap L}e_{k+1}(I,(X\cap L,x)).\tgs1.2$$
Note that the moving part cannot computed locally at points of $X\cap L,$
because a component $C$ of $D$ that gives rise to a moving component of the
Segre cycle maps to a subset of $X$ of codimension at most $k.$ So, the map
$$C\cap b\inv(L)\to X\cap L$$ has fibres of dimension at most $k-1.$ So, they
don't contribute to $e_{k+1}(I,(X\cap L,x))$ at any point of $x$ of
$X\cap L.$

\medskip
As the moving part of the Segre cycles of codimension $k$ come
from components of the exceptional divisor that map to subsets of smaller
codimension, controlling these components should be enough to control the
moving part of the Segre cycle of codimension $k.$ This is in fact true as the
next lemma shows.

As we will see in Proposition (\mink.3), the two assumptions in the following
lemma are very closely related.

\lem2 Assume that the integral closures of $I_1$ and $I_2$ coincide outside
a subset of codimension $k$ and
$e_i(I_1)=e_i^{i-1,1}(I_1,I_2)=e_i^{i-2,2}(I_1,I_2)$ for
$i=2,\dots,k-1.$ Then, we have
$$e_k(I_1)_m=e_k^{k-1,1}(I_1,I_2)_m=e_k^{k-2,2}(I_1,I_2)_m.$$

\pf Consider the family $$Y=P_{k-1}^{k-1,1}(I_1,tI_2;X\times\C)\to\C,$$ where
$t$ denotes the coordinate on $\C.$ Then, by \cite{GG, Proof of Corollary
(4.5)}, the assumptions imply that the fibre of $Y$ over 0 equals
$P_{k-1}(I_1,X).$

Now, consider a generic plane $L$ of codimension $n-k$ off 0 and the germ
$\tilde Y$ of $Y\cap (L\times\C)$ along the moving part of $Y\cap
(L\times\C)\cap V(I_1+tI_2),$ i.e. the part of this intersection that varies
with the tuple which is used to construct the polar variety. Now, by
assumption, the ideal $I_1^{[k-1]}+I_2^{[1]}$ generates a reduction of
$I_1|L$ and $I_2|L.$ Furthermore, by construction of $\tilde Y,$ we have
$$V(I_1^{[k-1]}+I_2^{[1]}|\tilde Y)=V(I_2^{[1]}|\tilde Y).$$
Hence, the degree of
$V(g_1|\tilde Y(0))$ equals
$e_k(I_1)_m,$ and, for non--zero $t,$ the degree of $V(g_1|\tilde Y(t))$
equals $e_k^{k-1,1}(I_1,I_2)_m.$  Hence, `conservation of number'
yields the first equality; see \cite{F, Prop. 10.2, p. 180}. The second
equality follows from an
analogous argument applied to
$$Y=P_{k-1}^{k-2,2}(I_1^{[k-1]}+I_1^{[1]},tI_1^{[2]};X\times\C)\to\C.$$\enddemo

\prop3 In the Setup (\segre.1), let $k$ be an integer so that $2\leq k\leq n$
 and assume that the equalities $e_1(I_1)=e_1^{1,1}(I_1,I_2)=e_1(I_2)$ hold,
and for
$j=2,\dots,k-1$ we have
$$\align
e_j(I_1)&=e_j^{j-1,1}(I_1,I_2)=e_j^{j-2,2}(I_1,I_2)\qquad\text{and}\\
e_j^{2,j-2}(I_1,I_2)&=e_j^{1,j-1}(I_1,I_2)=e_j(I_2).\endalign$$ Then, the
integral closure of
$I_1$ and $I_2$ coincide outside a subset of codimension $k.$ Also, the
following inequalities of mixed Segre numbers of codimension
$k$ obtain
$$\align e_k^{k-1,1}(I_1,I_2)^2&\leq
e_k(I_1)e_k^{k-2,2}(I_1,I_2)\qquad\text{and}\\
e_k^{1,k-1}(I_1,I_2)^2&\leq
e_k^{2,k-2}(I_1,I_2)e_k(I_2).\endalign$$

\pf We will proof the assertion by induction on $k.$ First, for any $k,$ the
proof of Corollary (\surf.6) carries over to the general case and shows that
the first assumed equality implies that $\bar I_1$ and $\bar I_2$ coincide
outside a subset of codimension two.

Assume now that the assertion is proven for $k-1.$ Then, by the above lemma,
we have the equalities
$$e_k(I_1)_m=e_k^{k-1,1}(I_1,I_2)_m=e_k^{k-2,2}(I_1,I_2)_m.$$
Hence, it is enough to show the claim for the multigerm of the intersection
of $X$ with a generic $(n-k)$--codimensional plane $L$ off 0 along the fixed
parts of the Segre cycles of codimension $k$ in question. Now, if $\bold h$
is a $(k-1)$--tuple of generic linear combinations of generators of
$I_1^{[k-1]}+I_2^{[1]},$ the polar curve
$P_{k-1}^{k-1,1}(I_1,I_2)$ equals the closure of
$V(h_1,\dots,h_{k-1})-V(I_1+I_2).$ But, by assumption, the integral closures
of the two ideals coincide outside a set of dimension 0. Hence, the
underlying sets of
$V(I_1+I_2)$ and $V(I_1)$ are equal. Also, without
changing $V(h_1,\dots,h_{k-1}),$ we may assume that
$h_2,\dots,h_{k-1}$ are generic linear combinations of $\bold f^{[k-1]}$
only. Furthermore, after changing $\bold f,$ we can write $h_2,\dots,h_{k-1}$
as generic linear combinations of $\bold f^{[k-2]}.$ Hence, the above polar
curve equals
$P_1^{1,1}(I_1,I_2;P_{k-2}(I_1)).$ A similar argument shows
$$P_{k-1}^{k-2,2}(I_1,I_2)=P_1(I_2,P_{k-2}(I_2)).$$ Hence, the first
inequality follows from our earlier result (\surf.5) on surfaces. Also, the
same result implies that the integral closures of the restrictions of $I_1$
and $I_2$ to $P_{k-2}(I_1)$ coincide. This implies already that $I_2$ is
integrally dependent on $I_1$ by a result of Teissier \cite{T3, Prop. (2),
p.349}. (Teissier assumes in his proof that $I_1$ is of finite colength and
that $X$ is Cohen-Macaulay, but it is easy to see that these assumptions are
unnecessary.) Then, the claim follows by symmetry.\enddemo

\medskip
The following
Corollary follows immediately from the previous proposition.

\cor4 In the Setup (\segre.1), the ideals $I_1$ and
$I_2$ have the same integral closure if and only if the equalities
$e_1(I_1)=e_1^{1,1}(I_1,I_2)=e_1(I_2)$ hold, and
for
$j=2,\dots,n$ we have
$$\align
e_j(I_1)&=e_j^{j-1,1}(I_1,I_2)=e_j^{j-2,2}(I_1,I_2)\qquad\text{and}\\
e_j^{2,j-2}(I_1,I_2)&=e_j^{1,j-1}(I_1,I_2)=e_j(I_2).\endalign$$
\endgroup\enddemo

\medskip
\flushpar{\it Example.} Let $X=(\C^3,0)$ with coordinates $x,y,z.$ Consider
the ideals
$I_1=(z)$ and $I_2=(xz,yz,z^2)$ in $\OO_{X,0}.$ Clearly, the two ideals
coincide outside 0. Also, we have $e_1(I_1)=e_1(I_2)=e_1^{1,1}(I_1,I_2).$
However, the polar surface of $I_1$ is empty, while the polar surface of $I_2$
isn't. In particular, we have $e_2(I_1)=0$ and $e_2(I_2)=1.$ So, the Segre
numbers of codimension 2 already distinguish $I_1$ and $I_2,$ even though
they're only different in codimension 3. In fact, the 2--codimensional Segre
cycle
$\Lambda_2(I_2)$ only has a moving component which is caused by the behavior
of the ideal at 0.
\medskip
Risler and Teissier \cite{T1, Ch.1,\S 2,\ p.302} proved a
product formula expressing the multiplicity of two ideals of finite colength
as a linear combination of the mixed multiplicities of the two ideals with
binomial coefficients. As the top Segre number $e_n(I)$ has a
length theoretic interpretation (see \cite{GG, (3.7)} and \cite{KT, (3.5)}),
one can derive the following product formula.

\prop5 In the setup (\segre.1), the following product formula obtains.
$$e_n(I_1I_2)=\sum_{i=0}^ne_n^{i,n-i}(I_1,I_2).$$

\pf The length theoretic interpretation of the top Segre number $e_n(J)$ of
an ideal $J\sub \OO_{X,0}$ shows that it equals the multiplicity of the ideal
$J_k=J\OO_{{X,0}_k}$ induced by
$J$ in the $k$th infinitesimal neighborhood of 0 in $X,$ for sufficiently big
$k.$ Now, the ideal $J_k$ is of finite colength, as any ideal in an Artinian
ring. Hence, we can apply the result of Risler and Teissier. This yields the
formula.\enddemo

\cor6 Under the assumptions of Proposition (\sectno.3), we have the following
product formula:
$$e_k(I_1I_2)=\sum_{i=0}^ke_k^{i,k-i}(I_1,I_2).$$

\pf As we have seen in the proof of Proposition (\sectno.3), the
assumptions imply that the moving part of $e_k^{i,k-i}(I_1,I_2)$ is
independent of $i.$ Furthermore, the same argument shows the moving part
$e_k(I_1I_2)_m$ equals
$e_k(I_1^2)_m.$ Now, the underlying sets of the blowups of $X$ along $I_1$ and
along $I_1^2$ are equal. Only, the exceptional divisor of the second blowup is
twice the divisor of the first blowup, and similarly for the first Chern
classes of the canonical bundles. Hence, from the alternative description of
Segre cycles discussed in (\sectno.1), we get

$$e_k(I_1^2)=2^ke_k(I_1).$$ Also, as the divisors of the two blowups have the
same underlying set, the same relation holds for the moving parts of these
numbers:

$$e_k(I_1^2)_m=2^ke_k(I_1)_m.$$ Hence, the claimed product formula holds for
the moving parts.

For the fixed parts, we intersect $X$ with a generic plane $L$ of codimension
$n-k$ off 0. Then, it is enough to show the equality locally at points of
$X\cap L.$ This, in turn, follows from the above proposition. \enddemo

\cor7 Under the assumptions of Proposition (\sectno.3), we have the following
Minkowski--type formula:
$$e_k(I_1I_2)^{1/k}\leq e_k(I_1)^{1/k}+e_k(I_2)^{1/k}.$$

\pf The proof follows almost exactly the proof of Teissier \cite{T2,
pp.39--42}. Only, we need to change his induction slightly. Instead of
considering $V(f_1)\subset X$ we need to consider $P_1^{f_1}(I_1)\subset X.$
Then, by the same argument as in the proof of Proposition (\sectno.3), we have
$$e_k^{i,k-i}(I_1,I_2)=e_{k-1}^{i-1,k-i}(I_1,I_2;P_1^{f_1}(I_1)).$$ The claim
then follows from this modification of Teissier's proof.\enddemo

\medskip
Teissier \cite{T3,Thm. 1} also showed that for ideals of finite colength this
Minkowski--type inequality is an equality if, and only if, there exist positive
integers $a$ and $b$ so that the integral closures of
$I_1^a$ and $I_2^b$ coincide. It is not possible to prove a similar result in
general, rather one would get a similar result with different powers in
different codimensions.

On the other hand, if the two ideals define subsets of $X$ of the same
dimension, and the components of Segre cycles of one of them form chains,
then we can derive an analogous result, using the observation after Corollary
(\surf.4).

\thm8 Let $(X,0)$ be an equidimensional analytic germ of dimension $n$ and
$I_2,I_2$ ideals in its local ring at 0. Assume that $V(I_1)$ and $V(I_2)$ are
of both of codimension $k.$ Furthermore, assume

$$|\Lambda_k(I_1)|\supset|\Lambda_{k+1}(I_1)|\supset
\dots\supset|\Lambda_n(I_1)|.$$
\itemitem{(1)} If $I_1$ and $I_2$ are of codimension at least two, then the
equality
$$ e_k(I_1I_2)^{1/k}=e_k(I_1)^{1/k}+e_k(I_2)^{1/k}$$ holds for
$k=2,\dots,n$ if and only if there exist positive integers $a$ and $b$ so that
the integral closures of
$I_1^a$ and $I_2^b$ coincide.
\itemitem{(2)} If $I_1$ and $I_2$ are of codimension one, then the
equalities $e_1(I_1)=e_1(I_1,I_2)=e_1(I_2)$ and
$$e_k(I_1I_2)^{1/k}=e_k(I_1)^{1/k}+e_k(I_2)^{1/k}$$ hold for $k=2,\dots,n$
if and only if the integral closures of
$I_1$ and $I_2$ coincide.

\pf We will only show that the Minkowski--type equalities implies the claimed
equalities of integral closures. The proof of the other direction follows
easily from the product formula (\sectno.6).

In case (1), all components of the exceptional
divisor of the blowup of $X$ along $I_1$ map to subsets of $X$ of codimension
at least $k.$ It follows that the Segre cycle $\Lambda_k(I_1)$ has no moving
components, and analogously for $I_2.$ Hence, the Segre numbers of codimension
$k$ can be computed locally at points of the intersection $X\cap L,$ where $L$
is a generic plane of codimension $n-k$ off 0. Now, Teissier showed \cite{T3,
p.354} that the assumed Minkowski--type equality imply
$$a^ke_k^{k,0}(I_1,I_2)=a^{k-1}be_k^{k-1,1}(I_1,I_2)=
\dots=b^ke_k^{0,k}(I_1,I_2).$$ Furthermore, a similar result to Lemma
(\mixed.2) holds for the $1/k$--norm. Hence, by Corollary (\mixed.4), we can
conclude that analogous equalities hold locally at points of $X\cap L.$ Now,
at a point
$x$ of $X\cap L$ the ideals
$I_1|(X\cap L,x)$ and $I_2|(X\cap L,x)$ are of finite colength. So, Teissier's
arguments imply
$$I_1^a|X\cap L=I_2^b|X\cap L.$$ Now, replacing $I_1$ by $I_1^a$ and $I_2$ by
$I_2^b,$ we  may assume that the integral closure of $I_1$ and $I_2$ coincide
outside a subset of codimension $k+1.$

In case (2), the argument of Corollary (\surf.6) apply again and show that
the integral closure of $I_1$ and $I_2$ coincide outside a subset of
codimension 2.

The following arguments apply now to both cases. We proceed by induction on
$k.$ So, assume now that $\bar I_1$ and $\bar I_2$ coincide outside a subset
of codimension $k-1.$ Then, by the same arguments as in the proof of Lemma
(\sectno.2) and Proposition (\sectno.3), it is enough to show the claim for
multi surface--germs. Furthermore, using Lemma (\mixed.2) once more, it
suffices to consider normal surface germs. Furthermore, the assumptions
imply that
the reduction to the surface case yields ideals of codimension one. This
case has been discussed already in Corollary (\surf.4); see also the remark
after the Corollary.\enddemo

\medskip
If we assume that one of the ideals contains the other, we can drop the
`chain--assumption' and get the generalization of Rees' theorem \cite{GG,
(4.9)}.

\thm9 In the Setup (\segre.1) assume $I_1\sub I_2.$
Then, the integral closures of $I_1$ and $I_2$ coincide if and only if the
equality $e_k(I_1)=e_k(I_2)$ holds for $k=1,\dots,n.$

\pf If the integral closures coincide, then, as the Segre numbers of an ideal
only depend on its integral closure, the desired equalities hold. So assume
now that the equalities of Segre numbers hold.

First, by the upper
semi--continuity of multiplicities applied to the ideals induced by
$I_1$ and
$I_2$ at points of
$X\cap L$ where
$L$ is a generic plane of codimension $n-1$ off 0, we have
$e_1(I_2)=e_1(I_1+I_2)\leq e_1(I_1).$ Hence, the assumptions imply that the
three numbers are in fact equal. Now, the remainder of the proof follows
exactly Teisier's original proof \cite{T3, p.355} for ideals of finite
colength: He shows that the assumptions imply that the Minkowski--type
equalities of Theorem (\sectno.8) hold. Then, as in the proof of the above
theorem, we show, by induction on the codimension $k,$ that the two ideals
have the same integral closure outside a subset of codimension $k.$ For the
induction step, again we can reduce to the case of surfaces. If the ideals we
end up with are of codimension one, then the result follows from
Corollary (\surf.4). And if they are of finite colength, we can apply
Teissier's original result.  This finishes the proof.\enddemo

\art10 An application in equisingularity theory

Consider two analytic function germs $f_0,f_1:(\C^{n+1},0)\to(\C,0)$ that
define reduced hypersurfaces $X(0)$ and $X(1).$ We view them as the fibres
over 0 and 1 of the family of hypersurfaces $X\sub(\C^{n+1},0)\times\C$
defined by $F(z,t)=f_0(z)+t(f_1(z)-f_0(z)),$ where $t$ is a coordinate of the
parameter axis $\C.$ We denote the singular locus of $X$ by $S(X).$ Let $m$ be
the maximal ideal in the local ring $\OO_{n+1}$ of $\C^{n+1}$ at 0.

In the theory of contact equivalence one defines the tangent space to the
contact equivalence class of a function $f$ on $(\C^{n+1},0)$ to be the ideal
$$T\Cal K(f)=mJ(f)+(f)\OO_{n+1}.$$ Here, $J(f)$ is the Jacobian ideal
generated by the partial derivatives of $f.$
Gaffney \cite{G1, Prop. 3.3} showed that the integral closure of this
tangent space controls some aspects of Whitney equisingularity of the family
$X:$

\smallskip
{\it The pair $(X-S(X),0\times\C)$ satisfies the Whitney conditions at
$(0,0)$ and $(0,1)$ if the integral closures of $T\Cal K(f_0)$ and $T\Cal
K(f_1)$ coincide.}

\smallskip
Using this result and our result on mixed Segre numbers we can give a
criterion for $X(0)$ and $X(1)$ to be members of a Whitney equisingular
family of hypersurfaces.

\prop11 In the above setup (\sectno.10), assume that for $k=2,\dots,n+1$ the
equalities
$$\align
e_k(T\Cal K(f_0))&=e_k^{k-1,1}(T\Cal K(f_0),T\Cal
K(f_1))=e_k^{k-2,2}(T\Cal K(f_0),T\Cal
K(f_1))\qquad\text{and}\\
e_k(T\Cal K(f_1))&=e_k^{k-1,1}(T\Cal K(f_1),T\Cal
K(f_0))=e_k^{k-2,2}(T\Cal K(f_1),T\Cal
K(f_0))\endalign$$ obtain. Then, there is
an Zariski--open subset $U$ of $\C$ containing the points 0,1 so that the
smooth part of
$X(U)=X\cap (U\times(\C^{n+1},0))$ is Whitney regular along $0\times U.$

\pf Gaffney's result together with Corollary (\mink.4) shows that the smooth
part of $X$ is Whitney regular along the parameter axis at the points $(0,0)$
and $(0,1).$ Also, as the Whitney--conditions hold generically, the existence
of $U$ as in the claim follows. This finishes the proof.\enddemo

\Refs
\atref{Fu}
\by William Fulton
\book Intersection theory
\publ Springer-Verlag
\yr 1984
\bookinfo Ergeb. Math.,3. Folge, 2. Band
\endref

\atref{G1}
\by Terence Gaffney
\paper Integral Closure of Modules and Whitney Equisingularity
\jour Invent. Math. \vol 107 \yr 1992 \pages 301--322
\endref

\atref{GG}
\by Terence Gaffney and Robert Gassler
\paper Numerical Invariants of Hypersurface SIngularities
\jour alg--geom e--print \vol 9611002\yr 1996
\endref

\ref \key KT
\by Steven L. Kleiman and Anders Thorup
\paper Mixed Buchsbaum--Rim multiplicity
\jour Kopenhagen Univ. Preprint \vol131\yr1994
\endref

\atref{T1}
\by Bernard Teissier
\paper Cycles \'evanescents, sections planes et conditions de Whitney
\jour Ast\'er\-isque \vol 7-8 \yr 1973 \pages 285--362
\endref

\atref{T2}
\by Bernard Teissier
\paper Sur une in\'egalit\'e \`a la Minkowski pour les multipliciti\'es
\jour Annals of Math.
\vol 106
\yr 1978 \pages 40--44
\endref

\atref{T3}
\by Bernard Teissier
\paper On a Minkowski--type inequality for Multiplicities -- II
\inbook C. P. Ramanujam -- A tribute \bookinfo Studies in Math. 8 \yr 1978
\pages 347--361
\publ Tata Institute
\endref

\endRefs
\enddocument